\documentclass[hidelinks,onefignum,onetabnum]{siamart251104}

\usepackage[dvipsnames]{xcolor}
\usepackage{hyperref}
\usepackage{cleveref}
\usepackage{url}
\usepackage{subcaption}
\usepackage{amsmath, amssymb, amsfonts}
\usepackage{graphicx}
\usepackage{comment}
\usepackage{booktabs}
\usepackage{multirow}
\usepackage{lipsum}
\usepackage{epstopdf}
\usepackage{algorithmic}
\usepackage{caption}
\usepackage{tikz}
\ifpdf
  \DeclareGraphicsExtensions{.eps,.pdf,.png,.jpg}
\else
  \DeclareGraphicsExtensions{.eps}
\fi

\newsiamremark{remark}{Remark}
\newsiamremark{hypothesis}{Hypothesis}
\crefname{hypothesis}{Hypothesis}{Hypotheses}
\newsiamthm{claim}{Claim}
\newsiamthm{assumption}{Assumption}
\crefname{assumption}{Assumption}{Assumptions}
\crefname{proposition}{Proposition}{Propositions}
\newsiamremark{fact}{Fact}
\crefname{fact}{Fact}{Facts}

\headers{Do DeepONet Trunk-Basis Preconditioners Need Training?}{J.~Oh, J.~Darbon and G.~Karniadakis}

\title{Do DeepONet Trunk-Basis Preconditioners Need Training?\thanks{Submitted to the editors \today.}
\funding{This work was supported by the U.S. Department of Energy (DOE), Office of Science, Office of Advanced Scientific Computing Research (ASCR), under Award Number DE-SC0023191 (SEA-CROGS: Scalable, Efficient and Accelerated Causal Reasoning Operators, Graphs and Spikes for Earth and Embedded Systems) and through the Scientific Discovery through Advanced Computing (SciDAC) Institute ``LEADS: LEarning-Accelerated Domain Science,'' under Subcontract \#831126 of Contract Number DE-AC05-76RL01830.}}

\author{Jaemin Oh\thanks{Division of Applied Mathematics, Brown University, Providence, RI 02912
  (\email{jaeminoh.math@gmail.com}, \url{https://jaeminoh.github.io}).}
\and J\'er\^ome Darbon\footnotemark[2]
\and George Em Karniadakis\footnotemark[2]
}

\usepackage{amsopn}

\ifpdf
\hypersetup{
  pdftitle={Do DeepONet Trunk-Basis Preconditioners Need Training?},
  pdfauthor={Jaemin Oh}
}
\fi

\begin{document}
\maketitle

\begin{abstract}
Trunk-basis (TB) preconditioners build the coarse space of a two-level Schwarz preconditioner from the trunk functions of a trained deep operator network (DeepONet), restricted to subdomains. Their effectiveness has been attributed to spectral bias since networks learn low frequencies first, complementing a smoother that removes the high-frequency error. Spectral bias, however, is a property of training, whereas a randomly initialized trunk with a smooth activation is already smooth. We therefore ask whether TB preconditioners need training. In one dimension, we prove that on each subdomain the span of untrained trunk functions is close to the space of low-degree polynomials under suitable nondegeneracy assumptions. Numerically, for the Poisson equation and three-dimensional linear elasticity with constant coefficients, untrained TB and piecewise polynomial coarse spaces of equal dimension yield close conjugate gradient iteration counts, while training with the standard \(\ell^2\) loss sometimes increases them. Training can improve iteration counts when high-contrast coefficients make the low-energy modes non-smooth, and if the loss measures the error in the energy norm. For the Poisson equation with contrast $10^5$, energy-trained TB coarse spaces need fewer iterations than piecewise polynomials of equal dimension, whereas \(\ell^2\)-trained ones do not. For a topology-optimized elastic cantilever with contrast $10^6$, the energy-trained TB space enriched with the rigid-body modes gives the fewest iterations at larger coarse dimensions. This suggests that for neural operator-based hybrid solvers, a network's inductive bias and training objective should match the solver.
\end{abstract}



\tableofcontents

\section{Introduction}
Neural operators learn the solution operator of a partial differential equation (PDE) from data~\cite{lu2021deeponet,li2021fno}. Their limited accuracy and generalization, however, can hinder their use in scientific and engineering applications~\cite{mcgreivy2024weak}. Recent work has therefore explored using neural operators to accelerate classical numerical solvers. Examples include neural operator warm-starts for Newton solvers~\cite{oh2026spectrally} and preconditioners for Krylov subspace methods~\cite{kopanicakova2025deeponet}. In this work, we focus on the latter approach.

Preconditioning is essential for the efficient solution of large-scale linear systems~\cite{saad2003iterative}. A suitable preconditioner improves the spectral properties of the system and accelerates the convergence of iterative methods~\cite[Lecture~40]{trefethen1997numerical}. One neural operator-based construction is the trunk-basis (TB) preconditioner developed in \cite{kopanicakova2025deeponet}. This approach restricts the trunk functions of a trained deep operator network (DeepONet)~\cite{lu2021deeponet} to the subdomains of a domain decomposition and uses their local spans to construct a coarse space for a two-level preconditioner~\cite[Section~3.4]{xu1992iterative}.

The effectiveness of the TB approach has been attributed to the spectral bias of neural networks~\cite{kopanicakova2025deeponet,zhang2024blending}. Specifically, it describes a phenomenon that multilayer perceptrons (MLPs) tend to fit low-frequency components before high-frequency components~\cite{rahaman2019spectral}. In this interpretation, the trained trunk network captures low-frequency features that complement classical relaxation methods (or smoothers) which efficiently reduce high-frequency error~\cite[Section~5.5]{xu2017algebraic}. The coarse correction in the TB space and the smoother thus act on complementary components of the error.

However, randomly initialized MLPs with smooth activation functions are already biased toward low-frequency functions~\cite{lee2018deep}. Training primarily adapts the trunk network to capture higher-frequency components, for which classical smoothers are already effective. This observation motivates us to investigate the TB space in the absence of training. Across our numerical experiments, training with the standard \(\ell^2\)-based loss often leads to an increase in the number of conjugate gradient (CG) iterations. This provides empirical evidence that the effectiveness of TB preconditioners arises from the intrinsic smoothness of randomly initialized trunk networks, rather than from training-induced adaptations.

We next investigate which functions are effectively represented by untrained TB spaces. On the theoretical side, we show that the span of trunk functions restricted to each subdomain approximates a low-degree polynomial space in one dimension, under suitable nondegeneracy assumptions. Numerically, we find that untrained TB spaces and piecewise polynomial coarse spaces of the same dimension yield comparable CG iteration counts for both Poisson's equation and linear elasticity with constant coefficients. These low-degree polynomial spaces contain the well-known smooth modes of the corresponding homogeneous equations,\footnote{Constant modes for Poisson's equation and rigid-body modes for linear elasticity.} explaining their effectiveness.

In general, effective coarse spaces must capture low-energy modes. For problems with low-contrast coefficients, these modes are typically low-frequency, whereas this is not necessarily the case for high-contrast problems. For high-contrast problems, low-energy modes may exhibit coefficient-dependent structure that local low-degree polynomials and smooth untrained trunk spaces approximate poorly. Training can then become beneficial by learning this structure. However, our experiments show that training with a standard \(\ell^2\)-based loss does not improve performance. An \(\ell^2\)-based objective may leave approximation errors with substantial energy, limiting the effectiveness of coarse correction. This motivates training in the energy norm, which weights errors in high-energy directions more strongly. For the Poisson equation with a high-contrast coefficient, energy-trained TB spaces outperform piecewise polynomial spaces of the same dimension. For high-contrast linear elasticity, combining energy-based training with rigid-body mode enrichment yields the best performance among the methods considered at larger coarse dimensions. The use of an energy-based objective is also consistent with~\cite{rudikov2024fcgno}, which trains a neural operator in the energy norm to approximate solutions of residual equations.

Our work complements recent studies of the design of neural operator-based hybrid solvers. In the context of replacing coarse corrections with neural operators, the importance of network architectures and training strategies is examined in \cite{melchers2026can,wu2026deep}, while \cite{hu2025hybrid} provides a convergence theory. However, we focus on the TB approach that extracts coarse spaces from learned neural operators. Random trunk-based preconditioning can be connected to \cite{avron2017faster}, where random feature-based preconditioners are used for kernel ridge regression. 

The rest of this paper is organized as follows. We review the mathematical preliminaries and TB preconditioners in \Cref{sec:background}. In \Cref{sec:low_contrast}, we examine why untrained TB spaces can be effective for constant-coefficient problems. In \Cref{sec:high_contrast}, we show that the energy norm-based loss function is more beneficial than the standard \(\ell^2\) loss for capturing low-energy modes of problems with high-contrast coefficients. We discuss limitations and future research directions in \Cref{sec:conclusion}.

\section{Background}\label{sec:background}

\subsection{Model problem}
Let \(\Omega \subset \mathbb{R}^d\) be a bounded Lipschitz domain whose boundary is partitioned into a Dirichlet part \(\Gamma_D\) of positive measure and a Neumann part \(\Gamma_N\). We consider second-order elliptic systems in divergence form for \(u : \Omega \to \mathbb{R}^{d_u}\),
\begin{equation}\label{eq:elliptic}
  -\nabla \cdot \bigl(K\, \mathcal{C} \nabla u\bigr) = f
  \quad \text{in } \Omega,
  \qquad
  u = 0 \quad \text{on } \Gamma_D,
  \qquad
  \bigl(K\, \mathcal{C} \nabla u\bigr)n = g
  \quad \text{on } \Gamma_N,
\end{equation}
where \(n\) is the outward unit normal, \(K : \Omega \to [K_{\min},K_{\max}]\) is a scalar coefficient with \(K_{\min}>0\), and \(\mathcal{C}\) is a constant linear map on \(\mathbb{R}^{d_u \times d}\). The ratio \(K_{\max}/K_{\min}\) is the contrast of \(K\). We focus on two cases:
\begin{itemize}
  \item Poisson's equation: \(d_u=1\) and \(\mathcal{C}\nabla u=\nabla u\);
  \item linear elasticity: \(d_u=d\) and
  \[
    \mathcal{C}\nabla u
    = 2\mu\,\varepsilon(u)+\lambda\,(\nabla\cdot u)I,
    \qquad
    \varepsilon(u)=\tfrac12(\nabla u+\nabla u^\top),
  \]
  where \(\varepsilon(u)\) is the infinitesimal strain tensor and
  \(\mu>0\), \(\lambda\geq 0\) are the Lam\'e constants.
\end{itemize}
The weak formulation of \cref{eq:elliptic} seeks
\(u\in V=\{v\in[H^1(\Omega)]^{d_u}:v=0\text{ on }\Gamma_D\}\)
such that
\begin{equation}\label{eq:weak_form}
  a(u,v)
  := \int_\Omega K\,\mathcal{C}\nabla u:\nabla v\,dx
  = \int_\Omega f\cdot v\,dx
    + \int_{\Gamma_N}g\cdot v\,ds
  \qquad \forall v\in V.
\end{equation}
For both problems, \(a\) is symmetric and coercive on \(V\). A conforming finite element space \(V_h\subset V\), with nodal basis \(\{\phi_i\}_{i=1}^N\), therefore yields a sparse symmetric positive definite (SPD) linear system
\begin{equation}\label{eq:linear_system}
  Au=f,
  \qquad
  A\in\mathbb{R}^{N\times N},
  \qquad
  A_{ij}=a(\phi_j,\phi_i),
  \qquad
  u,f\in\mathbb{R}^N.
\end{equation}
With a slight abuse of notation, \(u\) and \(f\) also denote the coefficient vector of the discrete solution and the assembled load vector, respectively. Systems of this form also arise from implicit time discretizations of parabolic PDEs.

Many-query applications require repeated solutions of \cref{eq:linear_system} for different instances of \(A\) and \(f\). In topology optimization, for example, the material coefficients change as the design is updated. For large \(N\), even a single solve can be expensive, and repeated linear solves can dominate the total computational cost.

The conjugate gradient (CG) method~\cite{hestenes1952methods} is a standard choice for SPD systems. Its classical convergence bound depends on the square root of the condition number of \(A\), so ill-conditioning can substantially increase the number of iterations. Preconditioning improves convergence by replacing the original system with
\[
    M^{-1} Au = M^{-1}f, \qquad M\in \mathbb{R}^{N \times N},
\]
where the preconditioner \(M^{-1}\) is computationally cheap to apply and the preconditioned operator \(M^{-1}A\) has more favorable spectral properties. Multigrid and domain decomposition methods provide widely used preconditioners for \(A\)~\cite{xu1992iterative,xu2017algebraic}.

\subsection{Two-level preconditioners}
Given an iterate \(u^k\), the exact error \(u-u^k\) satisfies the residual equation
\begin{equation}\label{eq:residual_equation}
  A(u-u^k)=r^k,
  \qquad r^k=f-Au^k.
\end{equation}
A preconditioner approximates the solution of this equation to obtain a correction. Two-level methods combine a local smoother with a coarse correction~\cite{xu2017algebraic}. The smoother reduces error components that are locally easy to resolve, while the coarse space approximates those that the smoother reduces slowly. For standard constant-coefficient problems, these components are often described as high- and low-frequency errors, respectively.

Let \(Q\in\mathbb{R}^{N\times N_c}\) have full column rank, with \(N_c\ll N\), and let its range define the coarse space. The Galerkin coarse operator is
\[
  A_c=Q^\top A Q\in\mathbb{R}^{N_c\times N_c}.
\]
Solving
\[
  A_c e_c^k=Q^\top r^k
\]
gives the coarse correction
\[
  e_{\mathrm{coarse}}^k
  =Qe_c^k
  =QA_c^{-1}Q^\top r^k.
\]
An effective coarse space must therefore approximate the error components left by the smoother. In linear elasticity, local rigid-body translations and rotations provide important coarse-space features. More generally, the relevant modes depend on the operator and may not be known a priori.

In this work, we use a symmetric multiplicative combination of smoothing and coarse correction,
\begin{equation}\label{eq:two_level}
  I-M^{-1}A
  =(I-\omega BA)\,
   (I-QA_c^{-1}Q^\top A)\,
   (I-\omega BA),
\end{equation}
where \(B\) is the smoother, such as Jacobi, Gauss--Seidel, and additive Schwarz methods, and \(\omega>0\) is a damping parameter. The factors on the right correspond to pre-smoothing, coarse correction, and post-smoothing.

We next review neural operators and the construction of coarse spaces from DeepONet trunk functions.

\subsection{Neural operators}
A PDE defines a solution operator that maps input functions or parameters to its solution. For simplicity, we write this map as \(\mathcal{G}: f \mapsto u\), with \(f\) denoting the input data. Neural operators approximate such maps between function spaces, although their implementation and training use finite-dimensional representations of the inputs and outputs. Among the available architectures~\cite{kovachki2023neural,lee2025finite,zhou2026transolver}, we focus on DeepONets~\cite{lu2021deeponet}. While we consider the map \(f \mapsto u\) in this work, other solution operators, e.g., \(K \mapsto u\), can be considered.

\subsubsection*{Deep operator networks}
A DeepONet combines a branch network \(\beta\), which encodes the input \(f\), and a trunk network \(\tau:\mathbb{R}^d\to\mathbb{R}^p\), which takes a spatial query point \(x\). For function-valued inputs, the branch typically receives values at fixed sensors. The prediction has the form
\begin{equation}\label{eq:deeponet}
  u_c(x)
  \approx \mathcal{G}_\theta[f]_c(x)
  = \sum_{q=1}^{p}\beta_{q,c}[f]\,\tau_q(x),
  \qquad c=1,\dots,d_u.
\end{equation}
The \(p\) trunk functions are shared across the \(d_u\) output components, and the branch returns \(pd_u\) coefficients. The network parameters \(\theta\) are typically learned from input--output pairs
\(\{(f_i,u_i)\}_{i=1}^{N_{\mathrm{train}}}\)
by minimizing a data-fitting loss, such as the discrete mean squared error
\[
  \mathcal{L}(\theta)
  = \frac{1}{N_{\mathrm{train}}}
    \sum_{i=1}^{N_{\mathrm{train}}}
    \|u_i-\mathcal{G}_\theta[f_i]\|_2^2,
\]
where the solution and prediction are evaluated at the training nodes.

Both data generation and network optimization can be costly. Collecting a sufficiently rich training dataset~\cite{shin2024twostep} requires repeated PDE solves, while optimization can be time-consuming for large networks and datasets. Moreover, a trained network may not achieve the accuracy required of a classical numerical solver~\cite{mcgreivy2024weak}. These limitations motivate hybrid methods that use neural operators to accelerate an iterative solver while retaining its convergence. The TB preconditioner~\cite{kopanicakova2025deeponet,kopanicakova2025leveraging} is one such approach.

\subsection{Trunk-basis preconditioners}\label{sec:tb}
Let \(\{x_1,\dots,x_n\}\subset\mathbb{R}^d\) be the discretization nodes. Evaluating the trunk functions at these nodes gives
\[
  T=
  \begin{bmatrix}
    \tau_1(x_1) & \tau_2(x_1) & \dots & \tau_p(x_1)\\
    \tau_1(x_2) & \tau_2(x_2) & \dots & \tau_p(x_2)\\
    \vdots & \vdots & & \vdots\\
    \tau_1(x_n) & \tau_2(x_n) & \dots & \tau_p(x_n)
  \end{bmatrix}
  \in\mathbb{R}^{n\times p}.
\]
For each component \(c\), the nodal prediction satisfies
\begin{equation}\label{eq:deeponet_matrix}
  \begin{bmatrix}
    u_c(x_1)\\
    \vdots\\
    u_c(x_n)
  \end{bmatrix}
  \approx
  \begin{bmatrix}
    \mathcal{G}_\theta[f]_c(x_1)\\
    \vdots\\
    \mathcal{G}_\theta[f]_c(x_n)
  \end{bmatrix}
  =
  T
  \begin{bmatrix}
    \beta_{1,c}[f]\\
    \vdots\\
    \beta_{p,c}[f]
  \end{bmatrix}.
\end{equation}
Thus, every predicted component lies in the range of \(T\), which is trained to approximate features shared by the solution data.

The construction in~\cite{kopanicakova2025deeponet} forms a coarse basis by selecting \(m\) columns of \(T\) at random, restricting them to subdomains of a domain decomposition~\cite{toselli2004domain}, and orthonormalizing the restrictions locally. Without loss of generality, let \(\tau_1,\dots,\tau_m\) be the selected functions. For the block construction considered here, partition the nodal indices into \(S\) disjoint groups, and denote the nodes in group \(s\) by \(\{x_{i,s}\}_{i=1}^{n_s}\). Assuming that each local restriction has full column rank, a reduced QR factorization gives
\[
  T_s=
  \begin{bmatrix}
    \tau_1(x_{1,s}) & \dots & \tau_m(x_{1,s})\\
    \vdots & & \vdots\\
    \tau_1(x_{n_s,s}) & \dots & \tau_m(x_{n_s,s})
  \end{bmatrix}
  =Q_sR_s.
\]
After ordering the nodes by subdomain, the scalar TB coarse basis is
\begin{equation}\label{eq:tb_space}
  Q=
  \begin{bmatrix}
    Q_1 & 0 & \dots & 0\\
    0 & Q_2 & \dots & 0\\
    \vdots & \vdots & \ddots & \vdots\\
    0 & 0 & \dots & Q_S
  \end{bmatrix},
  \qquad Q_s=T_sR_s^{-1}.
\end{equation}
For vector-valued problems, the scalar basis is applied independently to each component. The resulting coarse space has dimension \(d_u mS\), provided that all local blocks have full rank.

\subsection{Spectral bias explanation}\label{sec:spectral_bias}
The TB construction uses the trunk functions as a learned approximation space, while the branch predicts the coefficients of individual solutions; see \cref{eq:deeponet_matrix}. A well-trained trunk may therefore capture features common to a family of PDE solutions. TB-based methods have been applied successfully to Poisson's equation, linear elasticity, high-wavenumber Helmholtz problems, and the heat equation~\cite{kopanicakova2025deeponet,kopanicakova2025leveraging,lee2025helmholtz}.

Their effectiveness has been attributed to spectral bias~\cite{rahaman2019spectral}: MLPs tend to learn low-frequency components before high-frequency components. Under this interpretation, the trunk captures low-frequency features that complement a classical smoother, which efficiently reduces high-frequency error.\footnote{Other examples exploiting spectral bias in optimization include \cite{hoyer2019neural,oh2026effects}.} This explanation is intuitive when geometric frequency predicts the behavior of the smoother. For heterogeneous problems, however, the relevant notion of smoothness is determined by the discrete operator \(A\).

\subsection{Algebraic smoothness}\label{sec:algebraic_smoothness}
Let \(B\) be an SPD matrix and consider the generalized eigenvalue problem
\[
  A v_i = \lambda_i B^{-1} v_i,
  \qquad
  (v_i, v_j)_{B^{-1}} = \delta_{ij},
  \qquad
  0 < \lambda_1 \leq \dots \leq \lambda_N,
\]
where \((a, b)_{B^{-1}} = b^\top B^{-1}a\). Errors dominated by the modes associated with the smallest eigenvalues have a small generalized Rayleigh quotient,
\[
  E_{A, B^{-1}}(e) := 
  \frac{e^\top Ae}{e^\top B^{-1} e}
  =
  \frac{\sum_i\lambda_i(e, v_i)_{B^{-1}}^2}
       {\sum_i (e, v_i)_{B^{-1}}^2}.
\]
We call a nonzero vector \(e\) algebraically smooth relative to \(B\) when its generalized Rayleigh quotient \(E_{A,B^{-1}}(e)\) is small.

Effective coarse spaces, including TB spaces, should accurately approximate algebraically smooth error components that the smoother reduces slowly. However, algebraic smoothness does not imply geometric smoothness, particularly when the coefficients vary strongly.

\begin{remark}
    Our numerical experiments use an additive Schwarz smoother. For illustration, \Cref{fig:low_energy_1d,fig:poisson1d_trunk} show ordinary eigenvectors of \(A\), corresponding to \(B=I\) in the generalized eigenproblem above. These illustrate coefficient-dependent low-energy structure, but may not coincide with the slowly damped modes of the Schwarz smoother.
\end{remark}

\subsection{Motivating example}
The spectral-bias explanation leaves two distinctions unresolved. First, spectral bias concerns the dynamics of training, whereas the smoothness of a randomly initialized trunk is a property of its architecture and initialization. Second, a coarse space should approximate errors that are algebraically smooth, which may not be geometrically smooth.

To illustrate the second distinction, consider the one-dimensional diffusion operator
\[
  u\mapsto-(K u')' \qquad \text{on }(0,1),
\]
with homogeneous Dirichlet boundary conditions. We compare the constant coefficient \(K = 1\) with a high-contrast coefficient that equals \(10^5\) on selected intervals and \(1\) elsewhere.

\begin{figure}[ht]
  \centering
  \includegraphics[width=\textwidth]{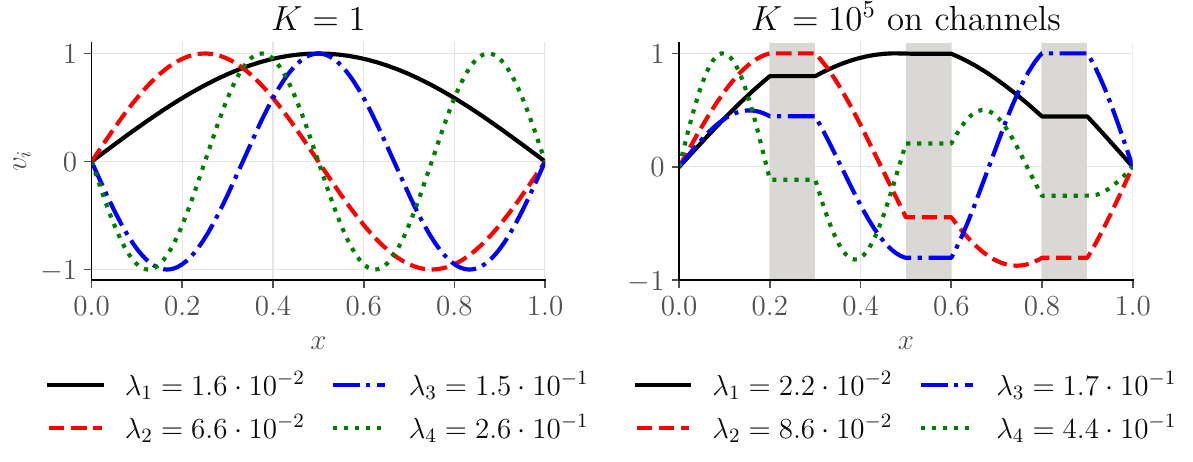}
  \caption{\textbf{Low-energy modes may not be geometrically smooth.}
  The first four eigenvectors of the discrete operator \(A\) for the one-dimensional problem \( - (K u')' = f\) on \((0,1)\), with homogeneous Dirichlet boundary conditions. Left: \(K = 1\). Right: \(K=10^5\) on the shaded intervals and \(1\) elsewhere. Coefficient contrast introduces localized features into the low-energy modes.}
  \label{fig:low_energy_1d}
\end{figure}

\begin{figure}[ht]
  \centering
  \includegraphics[width=\textwidth]{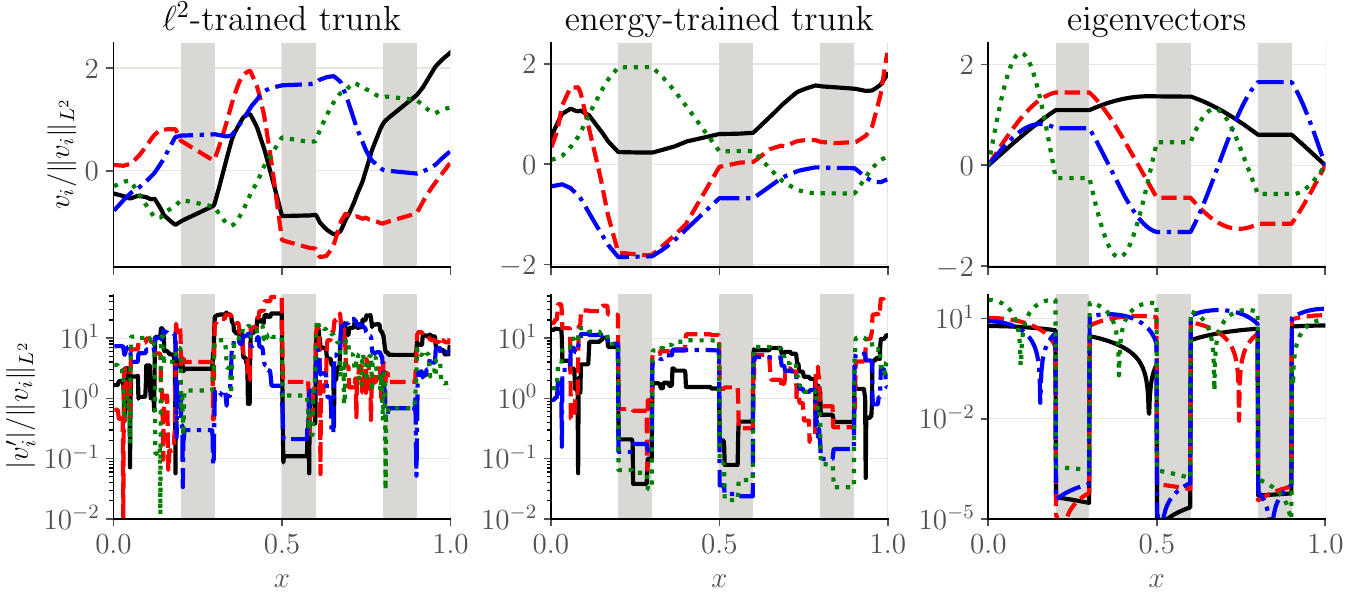}
  \caption{\textbf{Energy training captures coefficient-dependent structure.}
  The top panels show four vectors from the \(\ell^2\)-trained trunk, the energy-trained trunk, and the low-energy eigenspace; the bottom panels show the absolute values of their derivatives.
  The low-energy modes have nearly vanishing derivatives in the shaded high-conductivity regions.
  The energy-trained trunk vectors better reproduce this feature than the \(\ell^2\)-trained vectors.}
  \label{fig:poisson1d_trunk}
\end{figure}

\Cref{fig:low_energy_1d} shows the eigenvectors associated with the four smallest eigenvalues of the eigenvalue problem \(A v_i = \lambda_i v_i\). For constant coefficients, these modes are geometrically smooth. In the high-contrast case, they exhibit localized variations and changes in slope near coefficient interfaces. Low energy favors small derivatives within the high-conductivity regions, in which variations are strongly penalized. Effective coarse approximation must account for this coefficient-dependent structure. However, low-frequency features may not capture the localized structure of the low-energy modes.\footnote{Coefficient-aware domain decompositions or carefully placed biases in ReLU networks may help represent such features, but these approaches are beyond the scope of this work.}
Training can then be useful, provided that it learns the relevant structure without introducing excessive high-energy error. \Cref{fig:poisson1d_trunk} illustrates the role of the training objective: energy-trained trunk vectors better reproduce the small derivatives within high-conductivity regions than those trained with an \(\ell^2\)-based loss. These observations motivate the comparisons of untrained, \(\ell^2\)-trained, and energy-trained TB spaces in \Cref{sec:high_contrast}.

\section{Training can be unnecessary when low-energy modes are smooth}\label{sec:low_contrast}
For constant-coefficient problems, the error components that require coarse correction are often geometrically smooth. Familiar examples include low-frequency modes of Poisson's equation and local rigid-body motions in linear elasticity. For such problems, the smoothness of randomly initialized trunk functions may already provide the approximation properties needed for effective TB preconditioning. We first establish a local connection between TB spaces and polynomial spaces in one dimension, then examine its implications through numerical experiments in two and three dimensions.

\subsection{TB spaces approximate low-degree polynomials}\label{sec:tb_poly}
The TB construction in \cref{eq:tb_space} depends only on the local spans of the trunk functions. A randomly initialized MLP with a smooth activation produces smooth functions~\cite{lee2018deep}, whose restrictions to small subdomains admit accurate low-degree Taylor approximations. Under a suitable nondegeneracy condition, their span also approximates a local polynomial space. This provides a mechanism by which an untrained trunk can capture smooth coarse-space features.

\subsubsection*{The TB space in one dimension}
Let \(\Omega=(-1,1)\), with uniformly spaced nodes
\(-1=x_1<\dots<x_N=1\).
Let the trunk
\(\tau=(\tau_1,\dots,\tau_p):\mathbb{R}\to\mathbb{R}^p\)
be an MLP with a real-analytic activation, such as \(\tanh\).
The original TB construction selects \(m\leq p\) trunk functions at random~\cite{kopanicakova2025deeponet}. Without loss of generality, denote the selected functions by \(\tau_1,\dots,\tau_m\).

An alternative is to use the leading \(m\) left singular vectors of the nodal trunk matrix \(T\), thereby removing the additional variability introduced by random column selection. Provided that the corresponding singular values are nonzero, these vectors are nodal samples of linearly independent linear combinations of the original trunk functions. The resulting functions remain real analytic, so the argument below applies to this construction as well.\footnote{For \(N\geq p\), a dense thin SVD of \(T\in\mathbb{R}^{N\times p}\) requires \(\mathcal{O}(Np^2)\) operations.}
We therefore make the following assumption.

\begin{assumption}\label{ass:indep}
  The trunk functions \(\tau_1,\dots,\tau_m\) are linearly independent.
\end{assumption}

We partition the nodes into \(S\) contiguous groups \(I_s\), with
\(n_s:=|I_s|\geq\max\{m,2\}\).
Let \(\Omega_s=[c_s-h_s,c_s+h_s]\) be the smallest interval containing the nodes in \(I_s\), so that \(0<h_s\leq1\). The TB coarse space is
\[
  \mathcal{Q}=\bigoplus_{s=1}^{S}\mathcal{Q}_s,
  \qquad
  \mathcal{Q}_s
  :=\operatorname{span}
  \bigl\{(\tau_l(x_i))_{i\in I_s}:l=1,\dots,m\bigr\}
  \subset\mathbb{R}^{n_s}.
\]
For comparison, define the discrete polynomial spaces
\[
  \mathcal{P}_{s,r}
  :=\bigl\{(p(x_i))_{i\in I_s}:p\in\mathbb{P}_r\bigr\},
  \qquad
  \mathcal{P}_r:=\bigoplus_{s=1}^{S}\mathcal{P}_{s,r},
\]
where \(\mathbb{P}_r\) denotes the polynomials of degree at most \(r\).
Thus, \(\mathcal{P}_r\) is the discrete piecewise polynomial space with no continuity constraints between subdomains.

For subspaces \(X,Y\subset\mathbb{R}^{n_s}\), define the directed subspace gap
\begin{equation}\label{eq:delta}
  \delta(X,Y)
  :=\max_{x\in X,\,\|x\|_2=1}\|x-\Pi_Yx\|_2,
\end{equation}
where \(\Pi_Y\) is the Euclidean orthogonal projector onto \(Y\).
This quantity measures the largest relative best-approximation error for vectors in \(X\), and it vanishes if and only if \(X\subseteq Y\).
If \(U_X\) and \(U_Y\) have orthonormal columns spanning \(X\) and \(Y\), respectively, then
\[
  \delta(X,Y)
  =\bigl\|(I-U_YU_Y^\top)U_X\bigr\|_2,
\]
which can be computed by an SVD.

To state the local approximation result, define the Taylor coefficient matrix
\[
  G(x)_{j+1,l}
  :=\frac{\tau_l^{(j)}(x)}{j!},
  \qquad j=0,\dots,m-1,\quad l=1,\dots,m,
\]
and its singular set
\[
  Z:=\{x\in[-1,1]:\det G(x)=0\}.
\]
Under Assumption~\ref{ass:indep}, real analyticity implies that \(Z\) is finite; see \cref{lem:wronskian}.

\begin{proposition}[Local polynomial approximation by TB blocks]\label{prop:tb1d}
  Suppose that Assumption~\ref{ass:indep} holds and \(c_s\notin Z\).
  Then, for every \(r=0,\dots,m-1\),
  \begin{equation}\label{eq:tb1d}
    \delta\bigl(\mathcal{P}_{s,r},\mathcal{Q}_s\bigr)
    \leq C_s h_s^{m-r}.
  \end{equation}
  The constant \(C_s\) depends on the trunk functions and the conditioning of \(G(c_s)\), but not on \(n_s\).
  These constants are uniformly bounded when the subdomain centers remain a fixed positive distance from \(Z\).
\end{proposition}
\begin{proof*}
  See Appendix~\ref{app:proofs}.
\end{proof*}

\begin{figure}[ht]
  \centering
  \includegraphics[width=\textwidth]{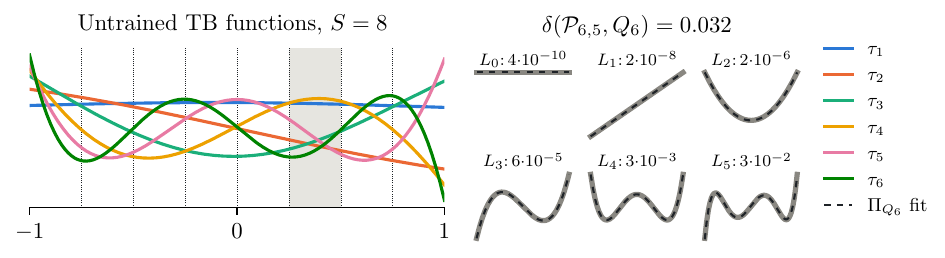}
  \caption{\textbf{Untrained trunk functions approximate low-degree polynomials locally.}
  Left: six leading left singular vectors of the nodal matrix of a randomly initialized trunk.
  Right: the first six Legendre polynomials on the highlighted subdomain, their approximations in the restricted trunk space, and the corresponding approximation errors.}
  \label{fig:trunk_1d}
\end{figure}

\Cref{prop:tb1d} shows that every vector of polynomial samples in
\(\mathcal{P}_{s,m-1}\) can be approximated in the TB block with relative Euclidean error at most \(C_sh_s\).
Lower-degree polynomials admit sharper bounds, with \(h_s^{m-r}\).
\Cref{fig:trunk_1d} illustrates this behavior for \(m=6\) and \(h_s=0.125\): the restricted trunk space accurately approximates the first six Legendre polynomials, with smaller errors for the lower-degree polynomials.

Because the blocks have disjoint supports, the local estimate extends directly to the global coarse space.

\begin{corollary}[Piecewise polynomial approximation by TB spaces]\label{cor:tb1d}
  Under the assumptions of \cref{prop:tb1d}, if \(c_s\notin Z\) for every \(s\), then
  \[
    \delta(\mathcal{P}_{m-1},\mathcal{Q})
    =\max_s\delta(\mathcal{P}_{s,m-1},\mathcal{Q}_s)
    \leq\max_s C_sh_s.
  \]
\end{corollary}
\begin{proof}
  Write \(v=\bigoplus_s v_s\in\mathcal{P}_{m-1}\), and set
  \(\delta_s=\delta(\mathcal{P}_{s,m-1},\mathcal{Q}_s)\).
  Since the orthogonal projection acts blockwise,
  \[
    \|v-\Pi_{\mathcal{Q}}v\|_2^2
    =\sum_s\|v_s-\Pi_{\mathcal{Q}_s}v_s\|_2^2
    \leq\bigl(\max_s\delta_s^2\bigr)\|v\|_2^2,
    \qquad  \because \sum_s \|v_s\|_2^2 = \|v\|_2^2.
  \]
  Equality in the resulting bound is attained by choosing a unit vector supported on a block that realizes the largest local gap. The final estimate follows from \cref{prop:tb1d}.
\end{proof}

The result is about Euclidean approximation of polynomial spaces in one dimension. Its implications for preconditioner performance in higher dimensions are examined numerically below.

\subsubsection*{Coarse spaces and solver configuration}
Motivated by this connection, we compare three coarse-space constructions:
\begin{itemize}
  \item \textbf{TB (\(\ell^2\) loss):} local restrictions of the leading \(m\) left singular vectors of the trunk matrix \(T\) after training with an \(\ell^2\)-based loss;
  \item \textbf{TB (untrained):} the same construction using \(T\) at random initialization;
  \item \textbf{\(\mathcal{P}_r\):} local polynomial spaces spanned by monomials of total degree at most \(r\).
\end{itemize}
All local bases are orthonormalized before constructing the coarse operator.
We consider \(r\in\{0,1,2,3\}\).
In \(d\) dimensions, the scalar polynomial space has
\(\binom{r+d}{d}\) basis functions per subdomain.
We therefore choose
\[
  m=\binom{r+d}{d}
\]
for the TB spaces, so that all three constructions have the same coarse-space dimension.
For vector-valued problems, each scalar basis is applied independently to every component.

We use CG with the symmetric multiplicative two-level Schwarz preconditioner in \cref{eq:two_level}.
The local operator \(B\) is an additive Schwarz method (ASM) with ILU(0) local solves and zero overlap on \(S=256\) subdomains. 
The coarse matrix \(A_c\) is factorized by LU~\cite{amestoy2001mumps}.

We consider two constant-coefficient instances of \cref{eq:elliptic}, both with \(K\equiv1\): Poisson's equation and linear elasticity.

\subsection{Poisson's equation}\label{sec:poisson2d_low}
On \(\Omega=(0,1)^2\), we consider
\begin{equation}\label{eq:poisson_2d}
  -\Delta u=f \quad \text{in }\Omega,
  \qquad
  u=0 \quad \text{on }\partial\Omega.
\end{equation}
This is the scalar case of \cref{eq:elliptic}, with \(K = 1\).
For the solver test, we use
\[
  f(x,y)=\sin(4\pi x)\sin(2\pi y)\sin(2\pi xy).
\]

\begin{figure}[ht]
  \centering
  \includegraphics[width=\textwidth]{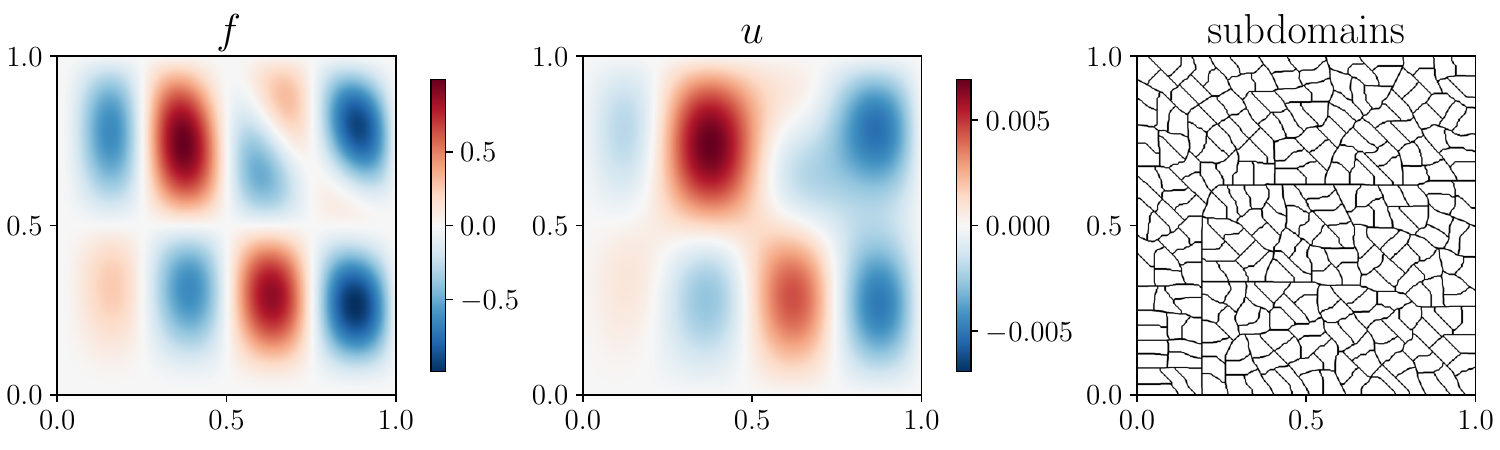}
  \caption{\textbf{Two-dimensional Poisson problem with \(K = 1\).} Left: test source term. Middle: corresponding numerical solution. Right: the partition into \(S=256\) subdomains used by the two-level Schwarz preconditioner.}
  \label{fig:poisson2d_low}
\end{figure}

We discretize the problem using Q1 elements on a uniform mesh \(\mathcal{T}^1\) with 49 squares per side and 2,500 degrees of freedom.
Five uniform refinements give \(\mathcal{T}^6\), with 2,461,761 degrees of freedom.
We generate the training data and train the DeepONet on \(\mathcal{T}^1\), then evaluate the preconditioners on \(\mathcal{T}^6\).
\Cref{fig:poisson2d_low} shows the test source term, its corresponding solution, and the partition into \(S=256\) subdomains.

To generate the dataset, we sample 5,000 Gaussian random fields \(f_i\) with covariance
\[
  k(x,y)=\exp\bigl(-|x-y|^2/(2\ell^2)\bigr),
  \qquad \ell=0.1,
\]
and solve \cref{eq:poisson_2d} on \(\mathcal{T}^1\) for the corresponding solutions \(u_i\). The last \(15\%\) of the samples are reserved for testing. We use the Adam optimizer~\cite{kingma2014adam}. Further training details are provided in Appendix~\ref{app:training}.

\begin{table}[t]
  \centering
  \caption{\textbf{CG iterations for the two-dimensional Poisson problem.} Piecewise polynomial, trained TB, and untrained TB coarse spaces are compared at equal dimensions, and with the same block sparsity pattern.}
  \begin{tabular}{lcccc}
    \toprule
    & \multicolumn{4}{c}{\(\dim Q\)} \\
    \cmidrule(lr){2-5}
    Coarse space & 256 & 768 & 1536 & 2560 \\
    \midrule
    \(\mathcal{P}_r\) & 150 & 79 & 51 & 35 \\
    TB (\(\ell^2\) loss) & 408 & 118 & 53 & 37 \\
    TB (untrained) & 149 & 79 & 51 & 35 \\
    \bottomrule
\end{tabular}
  \label{tab:poisson2d_low}
\end{table}

\Cref{tab:poisson2d_low} reports the number of CG iterations required to satisfy \(\|r^k\|_2\leq10^{-9}\). At every coarse dimension considered, the untrained TB space yields essentially the same iteration count as the piecewise polynomial space. The trained TB space requires more iterations at smaller dimensions and closely follows the untrained space only at larger dimensions, corresponding to \(\dim Q \in \{1536, 2560\}\). Thus, training does not improve preconditioner performance in this example.

\subsection{Linear elasticity}\label{sec:wrench}
We revisit the wrench example in~\cite[Section~4.2]{kopanicakova2025leveraging}.
The displacement \(u:\Omega\to\mathbb{R}^3\) satisfies
\begin{equation}\label{eq:elasticity}
  \begin{aligned}
    -\nabla\cdot\sigma(u)&=0
      &&\text{in }\Omega,\\
    u&=0
      &&\text{on }\Gamma_D,\\
    \sigma(u)n& = f
      &&\text{on }\Gamma_N,\\
    \sigma(u)n&=0
      &&\text{on }\partial\Omega\setminus(\Gamma_D\cup\Gamma_N),
  \end{aligned}
  \qquad
  \begin{aligned}
    \sigma(u)&=2\mu\,\varepsilon(u)
      +\lambda\,(\nabla\cdot u)I,\\
    \varepsilon(u)&=\tfrac12(\nabla u+\nabla u^\top).
  \end{aligned}
\end{equation}
Here, \(n\) is the outward unit normal.
We set Young's modulus to \(E=208\) and Poisson's ratio to \(\nu=0.3\), giving the Lam\'e constants \(\mu=80\) and \(\lambda=120\).

The wrench is clamped along the hexagonal hole \(\Gamma_D\) and subjected to the traction
\[
  f=(0,0,-g n_z)
\]
on the outer face \(\Gamma_N\) of its upper part, where \(n_z\) is the vertical component of \(n\).
The neural operator approximates the map \(g\mapsto u\).
Training inputs are drawn from \(g\sim U[-0.05,0.05]\), while all solver tests use \(g=0.1\).
Further training details are given in Appendix~\ref{app:training}.

\begin{figure}[ht]
  \centering
  \includegraphics[width=\textwidth]{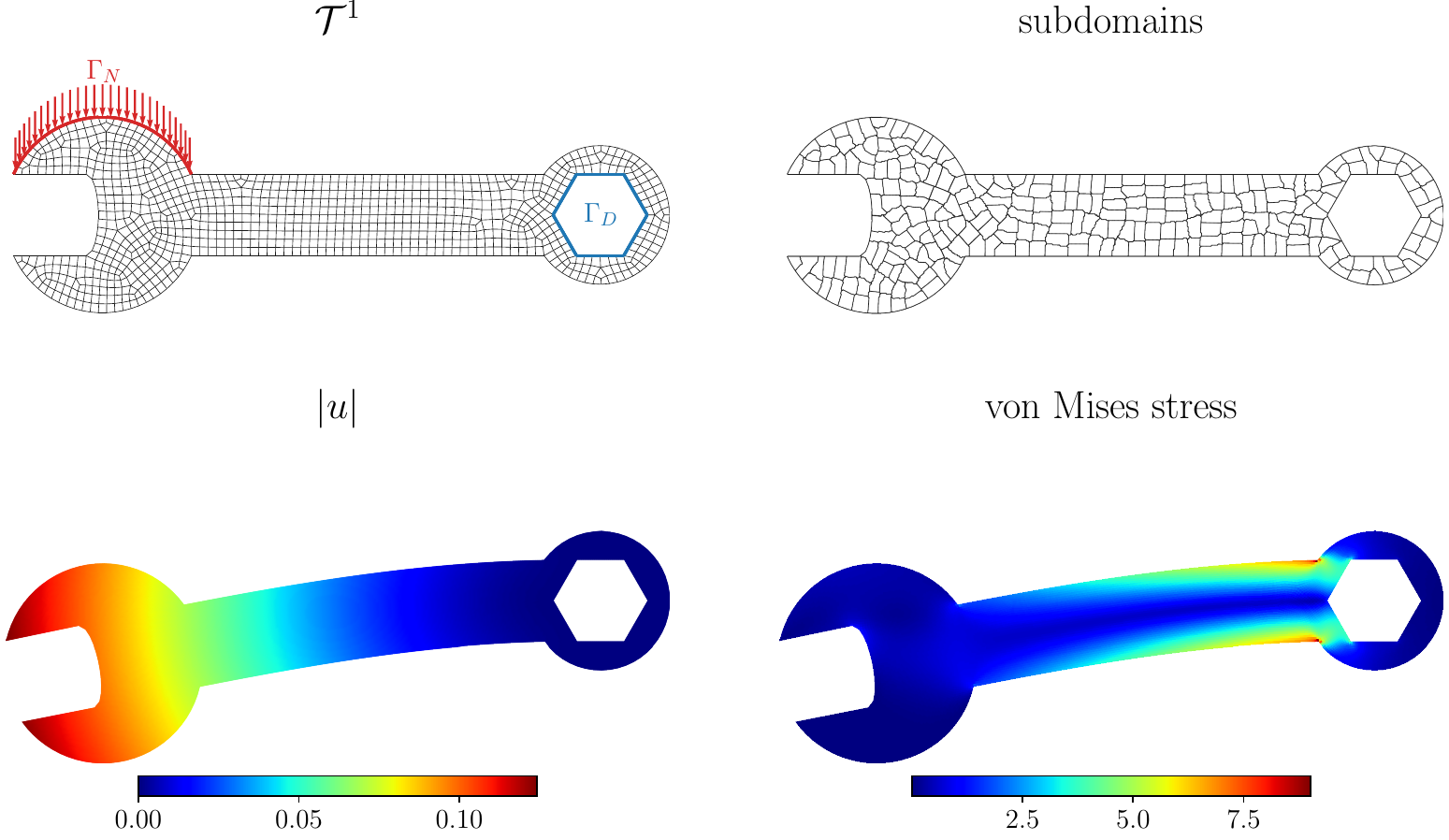}
  \caption{\textbf{Three-dimensional linear elasticity on a wrench.}
  Top: the mesh and boundary conditions, and the partition into \(S=256\) subdomains.
  Bottom: displacement magnitude and von Mises stress for \(g=0.1\).}
  \label{fig:wrench_overview}
\end{figure}

We discretize \cref{eq:elasticity} using Q1 elements.
The coarsest mesh \(\mathcal{T}^1\) contains 1,802 hexahedra and 9,216 degrees of freedom.
Three uniform refinements yield \(\mathcal{T}^4\), with 2,991,048 degrees of freedom, on which we evaluate the preconditioners.

\begin{table}[t]
  \centering
  \caption{\textbf{CG iterations for the wrench problem.}
  Piecewise polynomial, trained TB, and untrained TB coarse spaces are compared at equal dimensions and with the same block sparsity pattern.}
  \begin{tabular}{lcccc}
    \toprule
    & \multicolumn{4}{c}{\(\dim Q\)} \\
    \cmidrule(lr){2-5}
    Coarse space & 768 & 3072 & 7680 & 15360 \\
    \midrule
    \(\mathcal{P}_r\) & 614 & 60 & 29 & 23 \\
    TB (\(\ell^2\) loss) & 607 & 74 & 34 & 28 \\ 
    TB (untrained) & 616 & 61 & 31 & 28 \\
    \bottomrule
\end{tabular}
  \label{tab:wrench}
\end{table}

\Cref{tab:wrench} shows that the piecewise polynomial and untrained TB spaces yield similar iteration counts. Training increases the iteration count for \(\dim Q\in\{3072,7680\}\), while the difference between trained and untrained TB spaces decreases as the coarse dimension grows.

\begin{figure}[ht]
  \centering
  \includegraphics[width=\textwidth]{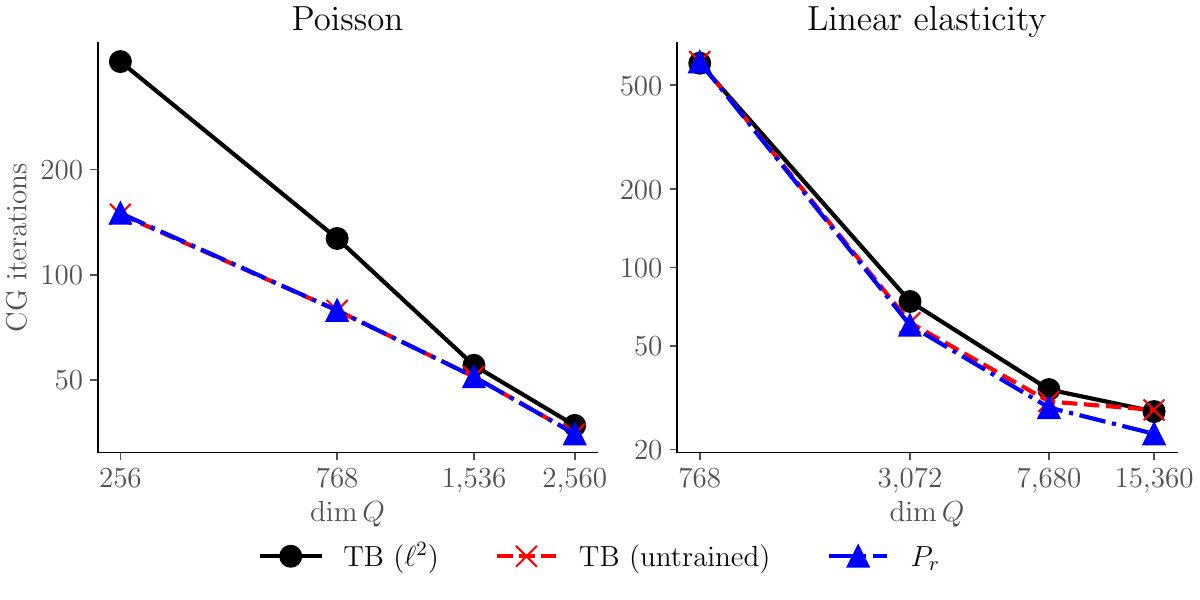}
  \caption{\textbf{Untrained TB and piecewise polynomial spaces yield similar iteration counts.}
  CG iterations as a function of the coarse-space dimension for the two constant-coefficient examples in \Cref{sec:low_contrast}.}
  \label{fig:why}
\end{figure}

Across two numerical examples (\cref{sec:poisson2d_low,sec:wrench}), we empirically observed excellent agreement between untrained TB and piecewise polynomial spaces as presented in \cref{fig:why}. This indicates that smooth untrained trunk functions can provide effective coarse spaces \emph{without} data collection or optimization, when low-energy modes are geometrically smooth.

The next section considers high-contrast problems, where low-energy modes have more complex spatial structure and training can offer benefits beyond those supplied by a smooth random initialization.

\section{The training objective matters when low-energy modes are not smooth}\label{sec:high_contrast}
Strong variations in the coefficient \(K\) can introduce localized structure into the low-energy modes, as illustrated in \cref{fig:low_energy_1d}. Piecewise polynomials and smooth random trunk functions may not approximate these modes well. Training offers a way to learn this coefficient-dependent structure from solution data.

To see why the data can contain the relevant modes, let
\(\{(\lambda_i,v_i)\}_{i=1}^N\) be the eigenpairs of \(A\), with Euclidean-orthonormal eigenvectors and
\(0<\lambda_1\leq\dots\leq\lambda_N\). Then
\begin{equation}\label{eq:sol_low_modes}
  u=A^{-1}f
  =\sum_{i=1}^{N}\frac{(f, v_i)}{\lambda_i}v_i.
\end{equation}
The inverse operator amplifies components of \(f\) associated with small eigenvalues. Thus, solution data can reveal low-energy modes, provided that the forcing distribution sufficiently excites them.\footnote{The forcing distribution is important: a mode absent from every training forcing is also absent from the corresponding solutions. A systematic study of training-data design is beyond the scope of this work.}
However, representing these modes accurately enough for coarse correction also depends on the training objective.

\subsection{Discrete energy loss}\label{sec:energy_loss}
An effective TB coarse space should approximate low-energy errors accurately in the energy norm. Even a perturbation with small Euclidean amplitude can substantially impair coarse correction if it lies in a high-energy direction. This motivates an energy-based training objective.

For an SPD matrix \(A\), define
\[
  \|x\|_A^2=x^\top Ax.
\]
For training pairs satisfying \(Au_k=f_k\), consider
\begin{equation}\label{eq:energy_loss}
    L_A(\theta)
  =\sum_k\|\mathcal{G}_\theta[f_k]-u_k\|_A^2.
\end{equation}
This objective is consistent with the energy-norm error minimized by CG. In the eigenbasis of \(A\),
\[
  \|\mathcal{G}_\theta[f_k]-u_k\|_A^2
  =\sum_i\lambda_i
    (\mathcal{G}_\theta[f_k]-u_k, v_i )^2.
\]
Unlike the unweighted \(\ell^2\) loss, which weights all eigenvector components equally, the energy loss assigns greater weight to errors associated with larger eigenvalues.

To connect this objective to coarse correction, let \(Q\) be a full-column-rank coarse basis and define the \(A\)-orthogonal projector
\[
  \Pi_Q^A=Q(Q^\top AQ)^{-1}Q^\top A.
\]
For a nonzero error \(e\), the fraction of its energy remaining after coarse correction is
\[
  \frac{\|(I-\Pi_Q^A)e\|_A^2}{\|e\|_A^2}
  =
  \frac{\min_c\|e-Qc\|_A^2}{\|e\|_A^2}.
\]
Effective correction therefore requires a small relative best-approximation error in the energy norm.

For illustration, consider eigenvectors \(v_L\) and \(v_H\), with eigenvalues \(\lambda_L\ll\lambda_H\), and the one-dimensional coarse space spanned by
\[
  q=v_L+\varepsilon v_H.
\]
Taking \(Q=[q]\) and \(e=v_L\) gives
\[
  \frac{\|(I-\Pi_Q^A)v_L\|_A^2}{\|v_L\|_A^2}
  =
  \frac{\varepsilon^2\lambda_H}
       {\lambda_L+\varepsilon^2\lambda_H}.
\]
Effective coarse correction requires \(\varepsilon^2\lambda_H\ll\lambda_L\); small perturbation magnitude \(|\varepsilon|\ll1\) alone is insufficient. This example explains how a small high-energy component can undermine an accurate \(\ell^2\) approximation of a low-energy mode.

\begin{remark}
    We note that \(v_L\) and \(v_H\) are usual eigenvectors, whereas we considered generalized eigenvectors in \Cref{sec:algebraic_smoothness}. This choice is for illustrative purposes.
\end{remark}

The energy loss also bounds the best-approximation error of the training solutions in the full trunk space. Let \(T\) denote the discrete trunk matrix, with its componentwise extension for vector-valued problems, and write
\(u_\theta[f]=T\beta_\theta[f]\). Then
\[
  \min_c\|u-Tc\|_A^2
  \leq\|u-T\beta_\theta[f]\|_A^2.
\]
Summing over the training samples shows that \(L_A\)~\cref{eq:energy_loss} bounds their aggregate best-approximation error in the trunk space. When the forcing distribution sufficiently excites the low-energy modes, energy-based training can therefore encourage their representation while suppressing high-energy contamination.

\subsection{Poisson's equation with high-contrast channels}\label{sec:poisson_contrast}
\begin{figure}[ht]
  \centering
  \includegraphics[width=\textwidth]{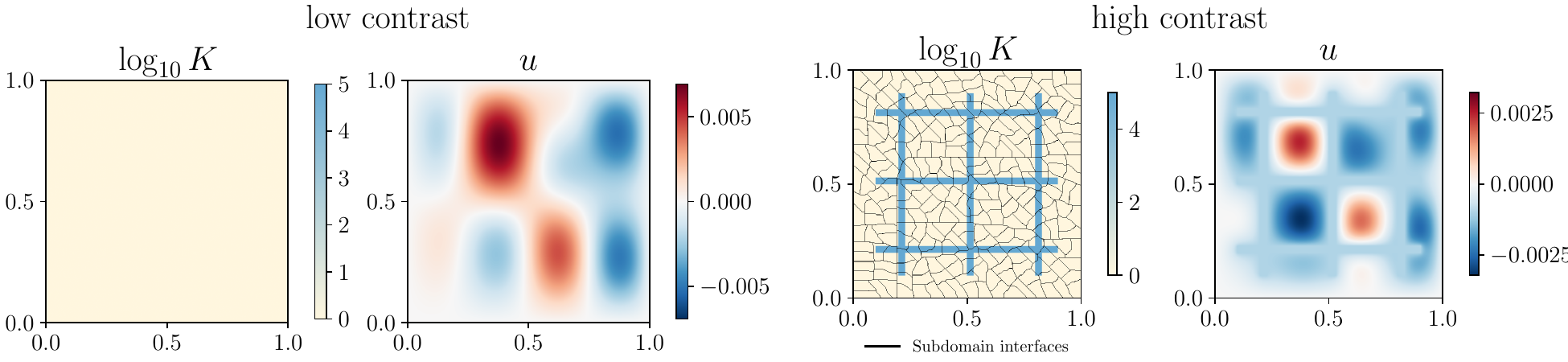}
  \caption{\textbf{Two-dimensional Poisson problems with constant and high-contrast coefficients.}
  Each pair of panels shows the coefficient field and a corresponding numerical solution.
  In the high-contrast case, \(K = 10^5\) on the channels and \(K = 1\) elsewhere. The domain decomposition is not aligned with the channels.}
  \label{fig:poisson_contrast}
\end{figure}

We consider
\[
- \nabla \cdot (K \nabla u) = f \quad \text{in }\Omega,
  \qquad
  u=0 \quad \text{on }\partial\Omega.
\]
The coefficient is \(K = 10^5\) on six channels of width \(0.03\) and \(K = 1\) elsewhere. The channels are aligned with the coordinate axes and centered at \(0.2\), \(0.5\), and \(0.8\) in each direction. \Cref{fig:poisson_contrast} shows the coefficient field and a corresponding solution. The test source term and the data generation procedure are the same as in \cref{sec:poisson2d_low}.

The \(\ell^2\)-trained network has a small relative \(\ell^2\) error and a large relative energy error~(\cref{tab:operator-errors}). The large energy error reflects variations of the prediction error along the high-conductivity channels, where gradients are strongly penalized. The energy-trained network has relative errors of \(0.46\) in both norms.

\begin{table}[t]
  \centering
  \caption{\textbf{CG iterations for constant- and high-contrast Poisson problems.} Piecewise polynomial and TB coarse spaces are compared at equal dimensions. Energy training reduces iteration counts for the high-contrast problem.}
  \label{tab:poisson_contrast}
  \resizebox{\textwidth}{!}{\begin{tabular}{lccccc}
    \toprule
    & \multicolumn{2}{c}{\(K = 1\)} & \multicolumn{3}{c}{\(K_{\max} / K_{\min} = 10^5\)} \\
    \cmidrule(lr){2-3} \cmidrule(lr){4-6}
    Coarse space & \(\dim Q = 768\) & \(\dim Q = 1536\) & \(\dim Q = 768\) & \(\dim Q = 1536\) & \(\dim Q = 2560\) \\
    \midrule
    \(\mathcal{P}_r\) & \bf 79 & \bf 51 & 198 & 172 & 137\\
    TB (untrained) & \bf 79 & \bf 51 & 198 & 172 & 137 \\
    TB (\(\ell^2\) loss) & 118 & 53 & 281 & 191 & 137 \\
    TB (energy loss) & 86 & 54 & \bf 164 & \bf 117 & \bf 107 \\
    \bottomrule
\end{tabular}
}
\end{table}

\Cref{tab:poisson_contrast} compares untrained, \(\ell^2\)-trained, and energy-trained TB spaces with \(\mathcal{P}_r\), using the solver configuration of \cref{sec:low_contrast}. Introducing the channels increases the iteration counts for every case. The untrained TB and polynomial spaces again yield similar counts: 198, 172, and 137 at the three coarse dimensions. The \(\ell^2\)-trained TB spaces require more iterations for \(\dim Q \in \{768, 1536\}\). Energy training reduces these counts to 164, 117, and 107, corresponding to reductions of approximately \(42\%\), \(39\%\), and \(22\%\) relative to the \(\ell^2\)-trained TB space.

For the constant-coefficient problem (\(K=1\)), energy training provides no improvement. The resulting spaces require 86 and 54 iterations, compared with 79 and 51 for \(\mathcal{P}_r\). In these experiments, the benefit of training therefore depends on both the coefficient-dependent structure of the low-energy modes and the choice of loss.

\subsection{Linear elasticity with high-contrast coefficients}\label{sec:cantilever}
High-contrast elasticity problems arise in heterogeneous materials and density-based topology optimization. We consider the latter, where a density field \(\rho\in[0,1]\) describes the material distribution, with \(\rho=0\) representing void and \(\rho=1\) representing solid material. A power-law interpolation involving \(\rho^p\), with penalization exponent \(p>1\), assigns much lower stiffness to low-density regions and produces a high-contrast coefficient field~\cite{bendsoe2004topology}. In our experiment, we use \(K = 10^{-6} + (1 - 10^{-6})\rho^3 \) exhibiting a contrast of \(10^6\).

\begin{figure}[ht]
  \centering
  \includegraphics[width=\textwidth]{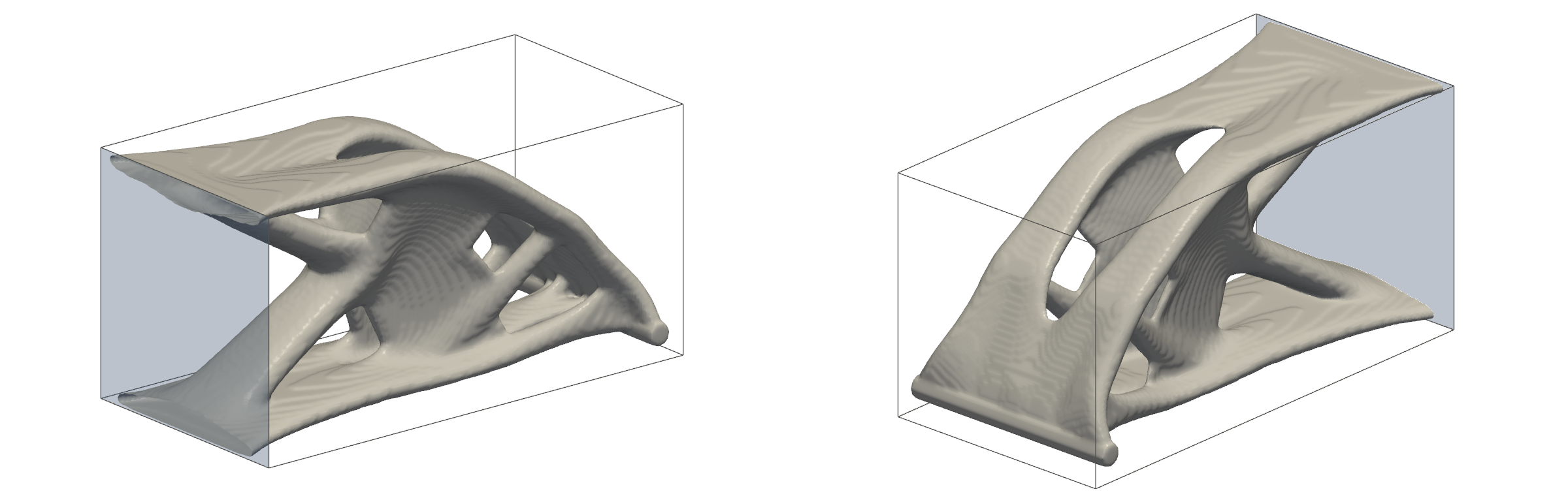}
  \caption{\textbf{High-contrast linear elasticity on a topology-optimized cantilever.}
  The surface shows the \(\rho=0.5\) isosurface of the optimized density field. The contrast is \(10^6\).}
  \label{fig:elasticity3d}
\end{figure}

We consider a three-dimensional cantilever beam of \([0, 2] \times [0, 1] \times [0, 1]\), fixed at one end and loaded at the other. A topology optimization algorithm~\cite{kim2025simple} produces the material distribution shown in \cref{fig:elasticity3d}, with a coefficient contrast of \(10^6\). The test problem is discretized using Q1 elements and has approximately 12.8 million degrees of freedom (\(256 \times 128 \times 128\) resolution). We compare untrained, \(\ell^2\)-trained, and energy-trained TB spaces on this fixed material distribution.

The \(\tanh\) trunk considered in our experiments did not adequately represent the complex spatial structure of this problem. We therefore use a SIREN trunk~\cite{sitzmann2020siren}, whose sinusoidal activations facilitate the representation of fine-scale features~\cite{oh2025spinnbgk,lee2025pinnp2d}. The frequency parameter \(\omega_0\) is set to \(10\).

We train on \(\mathcal{T}^1\), with 212,355 degrees of freedom (\(64 \times 32 \times 32\) resolution), using 500 forcing samples of the form
\[
  f=\left(\frac{D}{d_{\max}}\right)^{1/2}g,
  \quad
  D=\operatorname{diag}(A),
  \quad
  d_{\max}=\max_i D_{ii},
  \quad g\sim\mathrm{GRF}(\sigma=1, \ell=0.2).
\]
Both loss formulations use the SOAP optimizer~\cite{vyas2025soap}.

\begin{table}[t]
  \centering
  \caption{\textbf{CG iterations for the high-contrast cantilever problem.}
  The two-level Schwarz preconditioner uses \(S=256\) subdomains. TB spaces are compared with piecewise polynomial spaces, both with and without enrichment of the six rigid-body modes. The energy-trained TB spaces achieve \(\approx 25\%\) reductions in iteration counts compared to the \(\ell^2\) trained TB spaces for \(\dim Q = 7680\).}
  \label{tab:cantilever}
  \begin{tabular}{lcc}
    \toprule
    Coarse space & \(\dim Q = 3072\) & \(\dim Q = 7680\) \\
    \midrule
    NICO (\(\dim Q = 1536\)) & \multicolumn{2}{c}{327} \\
    \(P_1\), \(P_2\) & \bf 242 & 142 \\
    \midrule
    TB (untrained) & 329 & 144 \\
    TB (\(\ell^2\) loss) & 416 & 193 \\
    TB (energy loss) & 459 & 144 \\
    \midrule
    RBM + TB (untrained) & 260 & 141 \\
    RBM + TB (\(\ell^2\) loss) & 259 & 171 \\
    RBM + TB (energy loss) & 250 & \bf 129 \\
    \bottomrule
\end{tabular}

\end{table}

\Cref{tab:cantilever} reports iteration counts for the cantilever load.
At \(\dim Q=3072\), all three TB spaces require more iterations than the polynomial space: 329, 416, and 459 for the untrained, \(\ell^2\)-trained, and energy-trained spaces, respectively, compared with 242 for the polynomial baseline. The trained spaces also require more iterations than the coarse space with the six rigid-body modes (RBMs)~\cite{nicolaides1987deflation}, which takes 327 iterations with half the coarse dimension.

These results motivate explicit enrichment with the six local rigid-body modes. With this enrichment, the TB spaces require 260, 259, and 250 iterations at \(\dim Q=3072\), bringing their performance close to the polynomial baseline. At \(\dim Q=7680\), the enriched energy-trained space requires 129 iterations, compared with 142 for the polynomial space.
The enriched \(\ell^2\)-trained spaces remain less effective, requiring 259 and 171 iterations at the two dimensions.

\begin{table}[t]
  \centering
  \caption{\textbf{Additional polynomial and TB coarse-space comparisons for the cantilever problem.} CG iteration counts are reported for different numbers of subdomains and local basis sizes.}
  \begin{tabular}{llcccc}
    \toprule
    & & \multicolumn{4}{c}{\((r, m)\)} \\
    \cmidrule(lr){3-6}
    \(S\) & Coarse space & \((1, 2)\) & \((2,8)\) & \((3,18)\) & \((4,33)\) \\
    \midrule
    \multirow{2}{*}{256} & \(\mathcal{P}_r\) & 242 & 142 & 106 & 77 \\
        & RBM + TB (energy) & 250 & 129 & 91 & 70 \\
    \midrule
    \multirow{2}{*}{512} & \(\mathcal{P}_r\) & 179 & 108 & 82 & 60 \\
        & RBM + TB (energy) & 176 & 90 & 69 & 54 \\
    \bottomrule
\end{tabular}

  \label{tab:cantilever_poly_siren}
\end{table}

\Cref{tab:cantilever_poly_siren} systematically compares the piecewise polynomial space to the energy-trained TB space enriched with 6 rigid body modes. We vary numbers of subdomains \(S \in \{256, 512\}\),  polynomial degrees \(r \in \{1, 2, 3, 4\}\), and numbers of TB vectors \(m \in \{2, 8, 18, 33\}\) that yield the same coarse dimensions. Overall, the enriched, energy-trained TB space reduces iteration counts except for the case of \((S, r, m) = (256, 1, 2)\).

In this example, combining known mechanical modes with learned features gives the strongest result. The experiments highlight the complementary roles of coarse-space enrichment, dimension, and training objective.

\section{Conclusion}\label{sec:conclusion}
We investigated the mechanisms underlying the effectiveness of DeepONet trunk-basis (TB) preconditioners~\cite{kopanicakova2025deeponet} for SPD systems. Our analysis and experiments clarify the distinct roles of smoothness at initialization and subsequent training. In one dimension, we established that local TB spaces approximate low-degree polynomial spaces under a nondegeneracy condition. Experiments with constant-coefficient Poisson and elasticity problems showed that untrained TB and piecewise polynomial coarse spaces of equal dimension yield similar CG iteration counts, while training provides no consistent improvement. For high-contrast problems, low-energy modes can exhibit spatial structure that smooth random trunk functions do not readily capture. Training can then improve the coarse space, but its effectiveness depends on the loss. The energy-based objective controls high-energy perturbations in the TB space and yields more effective coarse spaces than the \(\ell^2\)-based objective in the examples considered. The high-contrast elasticity results further demonstrate the value of combining learned features with known rigid-body modes. Together, these findings support designing the network architecture, training objective, and coarse-space construction around the error components with which the classical solver struggles.

Several questions remain open. First, a quantitative theory is needed to connect the training objective and forcing distribution to the approximation properties of the localized, truncated trunk space and the resulting preconditioned spectrum. Second, extending the analysis to indefinite systems~\cite{murphy2000note}, including the Helmholtz equation~\cite{dolean2026achieving}, and studying its implications for time-dependent problems would broaden its applicability. Finally, practical assessments should account for state-of-the-art methods~\cite{vanek1996saamg,henson2002boomeramg,spillane2011robust}, data generation, training, coarse-space setup, and solve costs across repeated queries, together with robustness and parallel scalability.

\appendix
\section{Proof of \cref{prop:tb1d}}\label{app:proofs}

\begin{proof}\label{pf:prop_tb1d}
Fix a subdomain \(\Omega_s=[c_s-h_s,c_s+h_s]\), with \(h_s>0\), and relabel its nodes as \(x_1,\dots,x_{n_s}\). Define
\[
  \xi_i=\frac{x_i-c_s}{h_s}\in[-1,1].
\]
Because the physical nodes are uniformly spaced and \(\Omega_s\) is their smallest containing interval, the normalized nodes are uniformly spaced on \([-1,1]\).

Taylor expansion about \(c_s\) gives
\begin{equation}\label{eq:tb_1d_taylor}
  \tau_l(x_i)
  =\sum_{j=0}^{m-1}
    \frac{\tau_l^{(j)}(c_s)}{j!}
    h_s^j\xi_i^j
    +R_{il},
\end{equation}
where
\[
  |R_{il}|
  \leq\frac{M_m h_s^m}{m!},
  \qquad
  M_m:=\max_{1\leq l\leq m}
  \sup_{x\in[-1,1]}|\tau_l^{(m)}(x)|.
\]
Let
\[
  V^{(n_s)}
  :=
  \begin{bmatrix}
    1 & \xi_1 & \cdots & \xi_1^{m-1}\\
    1 & \xi_2 & \cdots & \xi_2^{m-1}\\
    \vdots & \vdots & & \vdots\\
    1 & \xi_{n_s} & \cdots & \xi_{n_s}^{m-1}
  \end{bmatrix},
  \qquad
  D_s:=\operatorname{diag}(1,h_s,\dots,h_s^{m-1}),
\]
and define
\[
  G(x)_{j+1,l}
  :=\frac{\tau_l^{(j)}(x)}{j!},
  \qquad
  j=0,\dots,m-1,\quad l=1,\dots,m.
\]
Writing \(G_s=G(c_s)\), the restricted trunk matrix satisfies
\[
  T_s=V^{(n_s)}D_sG_s+R,
  \qquad
  \|R\|_2\leq\|R\|_F
  \leq\sqrt{n_sm}\,\frac{M_mh_s^m}{m!}.
\]

Let \(v\in\mathcal{P}_{s,r}\), where \(0\leq r\leq m-1\).
Then \(v=V^{(n_s)}a\) for some \(a\in\mathbb{R}^m\) satisfying
\(a_{j+1}=0\) for \(j>r\).
Since \(c_s\notin Z\), the matrix \(G_s\) is invertible.
Choose
\[
  z=G_s^{-1}D_s^{-1}a.
\]
Then \(T_sz=v+Rz\in\mathcal{Q}_s\).
Because \(0<h_s\leq1\),
\[
  \|z\|_2
  \leq\|G_s^{-1}\|_2\|D_s^{-1}a\|_2
  \leq h_s^{-r}\|G_s^{-1}\|_2\|a\|_2.
\]
The best-approximation property of the orthogonal projector therefore gives
\begin{equation}\label{eq:bound_unnormalized}
  \begin{aligned}
    \|v-\Pi_{\mathcal{Q}_s}v\|_2
    &\leq\|v-T_sz\|_2\\
    &\leq\|R\|_2\|z\|_2\\
    &\leq
    \sqrt{n_sm}\,
    \frac{M_m}{m!}\,
    h_s^{m-r}\|G_s^{-1}\|_2\|a\|_2.
  \end{aligned}
\end{equation}

To normalize this estimate, define
\[
  \eta
  :=\inf_{k\geq\max\{m,2\}}
  \frac{\sigma_{\min}(V^{(k)})}{\sqrt{k}},
\]
where \(V^{(k)}\) is the Vandermonde matrix on \(k\) uniformly spaced nodes in \([-1,1]\).
We claim that \(\eta>0\).
Indeed, as \(k\to\infty\),
\[
  \frac{1}{k}(V^{(k)})^\top V^{(k)}
  \longrightarrow H,
  \qquad
  H_{j+1,l+1}
  =\frac12\int_{-1}^{1}x^{j+l}\,dx,
  \quad j,l=0,\dots,m-1.
\]
The matrix \(H\) is positive definite, since
\[
  a^\top Ha
  =\frac12\int_{-1}^{1}
    \left(\sum_{j=0}^{m-1}a_{j+1}x^j\right)^2dx
  >0
\]
for every nonzero \(a\).
Thus, the normalized smallest singular values are bounded away from zero for sufficiently large \(k\).
The remaining finitely many matrices have full column rank, which proves the claim.

Consequently,
\[
  \|v\|_2
  =\|V^{(n_s)}a\|_2
  \geq\eta\sqrt{n_s}\,\|a\|_2.
\]
Dividing \cref{eq:bound_unnormalized} by \(\|v\|_2\) and taking the supremum over nonzero \(v\in\mathcal{P}_{s,r}\) yields
\[
  \delta(\mathcal{P}_{s,r},\mathcal{Q}_s)
  \leq C_sh_s^{m-r},
  \qquad
  C_s
  :=\frac{\sqrt{m}\,M_m}{\eta\,m!}
    \|G(c_s)^{-1}\|_2.
\]
\end{proof}

For fixed trunk functions and \(m\), the dependence of \(C_s\) on the subdomain enters only through \(\|G(c_s)^{-1}\|_2\).
The following lemma controls this factor away from the singular set.

\begin{lemma}\label{lem:wronskian}
  Under Assumption~\ref{ass:indep} and the real-analyticity assumption on the trunk functions, the set
  \[
    Z:=\{x\in[-1,1]:\det G(x)=0\}
  \]
  is finite.
  For every fixed \(\rho>0\), the constants \(C_s\) are uniformly bounded over subdomains satisfying
  \(\operatorname{dist}(c_s,Z)\geq\rho\), independently of \(S\).
\end{lemma}

\begin{proof}[Proof of \cref{lem:wronskian}]\label{pf:wronskian}
The determinant of \(G\) is the Wronskian of
\(\tau_1,\dots,\tau_m\), multiplied by the nonzero constant
\(\prod_{j=0}^{m-1}(j!)^{-1}\).
For linearly independent real-analytic functions, the Wronskian is not identically zero.
Since the functions are analytic on a neighborhood of \([-1,1]\), the zeros of \(\det G\) are isolated and cannot accumulate in this compact interval.
Hence \(Z\) is finite.

For fixed \(\rho>0\), the set
\[
  K_\rho
  :=\{c\in[-1,1]:\operatorname{dist}(c,Z)\geq\rho\}
\]
is compact and contains no singular points of \(G\).
The map \(c\mapsto\|G(c)^{-1}\|_2\) is continuous on \(K_\rho\), so it is bounded there.
The expression for \(C_s\) then gives the asserted uniform bound.
If \(Z=\varnothing\), the same argument applies on all of \([-1,1]\).
\end{proof}

The lemma gives uniform local approximation estimates away from \(Z\).
If \(Z=\varnothing\), it also gives a constant \(C\), independent of the partition, such that
\[
  \delta(\mathcal{P}_{m-1},\mathcal{Q})
  \leq C\max_s h_s.
\]

\section{Details for DeepONet training}\label{app:training}
\begin{table}[t]
  \centering
  \caption{\textbf{DeepONet architectures and training hyperparameters.}
  Learning rates are constant for the Poisson and wrench problems; the cantilever learning rate decays from the listed initial value to zero following a cosine schedule~\cite{loshchilov2017sgdr}.}
  \label{tab:hyperparameters}
  \resizebox{\textwidth}{!}{
  \begin{tabular}{llll}
    \toprule
    & \Cref{sec:poisson2d_low,sec:poisson_contrast} & \Cref{sec:wrench} & \Cref{sec:cantilever}\\
    \midrule
    Trunk & \([2,256,256,128]\) & \([3,512,256,256]\) & \([3,512,512,256]\)\\
    Trunk activation & \(\tanh\) & \(\tanh\) & \(\sin(10x)\)\\
    Branch & \([2500,256,256,256,128]\) & \([1,256,256,256,768]\) & \([425,768]\)\\
    Branch activation & \(\tanh\) & \(\tanh\) & None (linear)\\
    Optimizer & Adam & Adam & SOAP\\
    Maximum epochs & 50,000 & 100,000 & 50,000\\
    Learning rate & \(10^{-3}\) & \(10^{-4}\) & \(10^{-4}\)\\
    \# of samples & 5,000 & 10,000 & 500 \\
    Batch size & 250 & 1,000 & 85\\
    \bottomrule
  \end{tabular}
  }
\end{table}
\Cref{tab:hyperparameters} summarizes the network architectures and training hyperparameters. The listed widths include the input and output layers. For each problem, the \(\ell^2\)-based and energy-based training runs use the same architecture and optimization hyperparameters. Training terminates after 10,000 epochs without improvement in the monitored loss or upon reaching the maximum epoch count.

\begin{table}[t]
    \centering
    \caption{\textbf{Mean relative errors in $\ell^2$ and energy norms.}}
    \label{tab:operator-errors}
    \begin{tabular}{llcccc}
        \toprule
        & & \multicolumn{2}{c}{$\ell^2$ error}
          & \multicolumn{2}{c}{energy error} \\
        \cmidrule(lr){3-4}\cmidrule(lr){5-6}
        Problem & Loss & Train & Test & Train & Test \\
        \midrule
        \multirow{2}{*}{\Cref{sec:poisson2d_low}}
        & \(\ell^2\)             & 0.0243  & 0.334   & 0.0882 & 0.357 \\
        & energy          & 0.0539  & 0.421   & 0.0998 & 0.508 \\
        \hline
        \Cref{sec:wrench} & \(\ell^2\) & 0.0165  & 0.0123  & 3.20   & 2.22  \\
        \hline
        \multirow{2}{*}{\Cref{sec:poisson_contrast}}
        & \(\ell^2\)             & 0.0384  & 0.324   & 25.3   & 25.7  \\
        & energy          & 0.0890  & 0.457   & 0.187  & 0.461 \\
        \hline
        \multirow{2}{*}{\Cref{sec:cantilever}}
            & \(\ell^2\)             & 0.0334  & 1.34    & 7.31   & 9.03  \\
            & energy            & 0.364   & 1.19    & 0.225  & 1.52  \\
        \bottomrule
    \end{tabular}
\end{table}
\Cref{tab:operator-errors} summarizes the final train and test relative errors in discrete \(\ell^2\) norm or energy norm.

\section{Seed variability for untrained TB coarse spaces}\label{app:seed}
Untrained TB spaces may be sensitive to random seeds for their initializations. To evaluate the sensitivity, we measure the iteration counts five times with the seed from \(0\) to \(4\). We use the same network architecture outlined in \cref{tab:hyperparameters}. \Cref{tab:seed} presents the results.

\begin{table}[t]
    \centering
    \caption{\textbf{Random seed variability for untrained TB preconditioning.}
    \(^*\)RBM + TB; all other rows use TB.}
    \begin{tabular}{lrcccccr}
    \toprule
    \multirow{2}{*}{Problem} & \multirow{2}{*}{\(\dim Q\)} & \multicolumn{5}{c}{Seed} & \multirow{2}{*}{Average} \\
    \cmidrule(lr){3-7}
     & & 0 & 1 & 2 & 3 & 4 & \\
    \midrule
    \multirow{4}{*}{\Cref{sec:poisson2d_low}} & 256 & 149 & 149 & 149 & 149 & 149 & 149.0 \\
     & 768 & 79 & 79 & 79 & 79 & 79 & 79.0 \\
     & 1536 & 51 & 51 & 51 & 51 & 51 & 51.0 \\
     & 2560 & 35 & 35 & 35 & 35 & 35 & 35.0 \\
    \midrule
    \multirow{4}{*}{\Cref{sec:wrench}} & 768 & 616 & 616 & 616 & 615 & 616 & 615.8 \\
     & 3072 & 61 & 62 & 62 & 61 & 61 & 61.4 \\
     & 7680 & 30 & 31 & 30 & 32 & 30 & 30.6 \\
     & 15360 & 28 & 28 & 28 & 29 & 29 & 28.4 \\
    \midrule
    \multirow{2}{*}{\Cref{sec:poisson_contrast}} & 768 & 198 & 198 & 198 & 198 & 198 & 198.0 \\
     & 1536 & 172 & 172 & 172 & 172 & 172 & 172.0 \\
    \midrule
    \multirow{4}{*}{\Cref{sec:cantilever}} & 3072 & 346 & 329 & 289 & 325 & 388 & 335.4 \\
     & 7680 & 143 & 144 & 143 & 145 & 141 & 143.2 \\
     & \(3072^*\) & 252 & 260 & 244 & 248 & 259 & 252.6 \\
     & \(7680^*\) & 142 & 141 & 142 & 140 & 140 & 141.0 \\
    \bottomrule
\end{tabular}

    \label{tab:seed}
\end{table}

\section{Computational environment}
The numerical experiments are implemented in Python. Finite element discretization and assembly use FEniCSx (DOLFINx)~\cite{baratta2023dolfinx}. Sparse linear algebra, Krylov solvers, and preconditioners use PETSc through its Python interface, petsc4py~\cite{petsc4py}. Linear solves run on CPUs, with MPI used for the distributed-memory cantilever calculations. The DeepONets are implemented using JAX~\cite{jax2018github,deepmind2020jax,kidger2021equinox} and NumPy~\cite{harris2020array}. Matplotlib~\cite{hunter2007matplotlib} is used to generate the figures.

\section*{Acknowledgement}
We thank \href{https://kopanicakova.github.io}{Alena Kopani{\v{c}}{\'a}kov{\'a}} and \href{https://sites.google.com/view/youngkyulee/home}{Youngkyu Lee} for sharing their code for numerical examples. Part of this research was conducted using the computational resources and services provided by the Center for Computation and Visualization at Brown University.

\bibliographystyle{modules/siamplain.bst}
\bibliography{references.bib}
\end{document}